\documentclass[12pt]{article}
\usepackage{amsmath}
\usepackage{amsopn}
\DeclareMathOperator{\tg}{tg}
\DeclareMathOperator{\tang}{tang}
\begin{document}
\title{A commentary on the continued fraction by which the illustrious La
Grange has expressed the binomial powers
\footnote{Delivered to the St.--Petersburg Academy March 20, 1780. Originally published as
\emph{
Commentatio in fractionem continuam, qua illustris La Grange potestates
binomiales expressit},
M\'emoires de l'acad\'emie des sciences de St.--Petersbourg \textbf{6} (1818),
3--11, and
republished in \emph{Leonhard Euler, Opera Omnia}, Series 1:
Opera mathematica,
Volume 16, Birkh\"auser, 1992. A copy of the original text is available
electronically at the Euler Archive, at http://www.eulerarchive.org. This paper
is E750 in the Enestr\"om index.}}
\author{Leonhard Euler\footnote{Date of translation: July 21, 2005.
Translated from the Latin
by Jordan Bell, 3rd year undergraduate in Honours Mathematics, School of Mathematics and Statistics, Carleton University,
Ottawa, Ontario, Canada.
Email: jbell3@connect.carleton.ca.
This translation was written
during an NSERC USRA supervised by Dr. B. Stevens.
}}
\date{}
\maketitle

I. This illustrious man has converted the Binomial power
$(1+x)^n$, by his most singular method of logarithmic differentials,
into
this continued fraction: 
\[
  (1+x)^n = 1 + \frac{nx}{1 
          + \frac{(1-n)x}{2 
          + \frac{(1+n)x}{3
          + \frac{(2-n)x}{2
          + \frac{(2+n)x}{5
          + \frac{(3-n)x}{2
          + \frac{(3+n)x}{7
          + \textrm{ etc.}}}}}}}},
\]
which expression celebrates the marvelous property that whenever 
the exponent $n$ is an integral number, either positive or negative,
it is halted and is reduced to a finite form.

II. Seeing that this continued fraction does not proceed by a uniform
law, but rather is interrupted,
we may bring it to a uniform law,
as it would be most desirable, if we should represent it in the following
way by parts:

\begin{center}
\begin{tabular}{rcl}
$(1+x)^n$&$=$&$1+\frac{nx}{A}$;\\
$A$&$=$&$1+\frac{(1-n)x}{2+\frac{(1+n)x}{B}}$;\\
$B$&$=$&$3+\frac{(2-n)x}{2+\frac{(2+n)x}{C}}$;\\
$C$&$=$&$5+\frac{(3-n)x}{2+\frac{(3+n)x}{D}}$;\\
$D$&$=$&$7+\frac{(4-n)x}{2+\frac{(4+n)x}{E}}$;\\
&etc.&
\end{tabular}
\end{center}

From here, we therefore will have by reducing:

\begin{center}
\begin{tabular}{lll}
$A$&$=1+\frac{(1-n)Bx}{2B+(1+n)x}$&$=1+\frac{(1-n)x}{2}-
\frac{(1-nn)xx\cdot 2}{2B+(1+n)x}$\\
&&$=1+\frac{(1-n)x}{2}+\frac{(nn-1)xx\cdot 4}{B+(\frac{1+n}{2})x}$.
\end{tabular}
\end{center}

In a similar way it will be:

\begin{center}
\begin{tabular}{lll}
$B$&$=3+\frac{(2-n)Cx}{2C+(2+n)x}$&$=3+\frac{(2-n)x}{2}-
\frac{(4-nn)xx\cdot 2}{2C+(2+n)x}$\\
&&$=3+\frac{(2-n)x}{2}+\frac{(nn-4)xx\cdot 4}{C+(\frac{2+n}{2})x}$.
\end{tabular}
\end{center}

In the very same way we shall have:

\begin{center}
\begin{tabular}{lll}
$C$&$=5+\frac{(3-n)Dx}{2D+(3+n)x}$&$=5+\frac{(3-n)x}{2}-
\frac{(9-nn)xx\cdot 2}{2D+(3+n)x}$\\
&&$=5+\frac{(3-n)x}{2}+\frac{(nn-9)xx\cdot 4}{D+(\frac{3+n}{2})x}$,
\end{tabular}
\end{center}

and so on.

III. Now if we were to substitute these values in succession in the place
of $A,B,C$, etc., the continued fraction will be induced to the following
form:
\[
  (1+x)^n = 1 + \frac{nx}{1 + \frac{(1-n)x}{2}
          + \frac{(nn-1)xx\cdot 4}{3(1+\frac{1}{2}x)
          + \frac{(nn-4)xx\cdot 4}{5(1+\frac{1}{2}x)
          + \frac{(nn-9)xx\cdot 4}{7(1+\frac{1}{2}x)
          + \frac{(nn-16)xx\cdot 4}{\textrm{ etc.}
          }}}}},
\]

IV. So that we can remove these partial fractions, we set $x=2y$,
so that this expression shall be obtained:
\[
  (1+2y)^n = 1 + \frac{2ny}{1 + (1-n)y
          + \frac{(nn-1)yy}{3(1+y)
          + \frac{(nn-4)yy}{5(1+y)
          + \frac{(nn-9)yy}{7(1+y) + \textrm{ etc.}}}}},
\]
which form can easily be transformed into this:
\[
\frac{2ny}{(1+2y)^n-1}=1+(1-n)y+\frac{(nn-1)yy}{3(1+y)+\frac{(nn-4)yy}{5(1+y)+\textrm{ etc.}}}.
\]
Then $ny$ is added to both sides, so that it will emerge
\[
\frac{ny(1+(1+2y)^n)}{(1+2y)^n-1}=1+y+\frac{(nn-1)yy}{3(1+y)+
\frac{(nn-4)yy}{5(1+y)+\textrm{ etc.}}};
\]
this expression is ordered enough that it may proceed to be
regulated.

V. We will now divide both sides by $1+y$, and the left
member will come out as: $\frac{ny}{1+y}\cdot \frac{(1+2y)^n+1}{(1+2y)^n-1}$.
From the right side moreover, each of the fractions on the top and bottom
should be divided by $1+y$, and this form will extend:
\[
  1 + \frac{(nn-1)yy\cdot (1+y)^2}{3
          + \frac{(nn-4)yy\cdot (1+y)^2}{5
          + \frac{(nn-9)yy\cdot (1+y)^2}{7
          + \frac{(nn-16)yy\cdot (1+y)^2}{9
          + \frac{(nn-25)yy\cdot (1+y)^2}{11
          + \textrm{ etc.}}}}}},
\]

VI. Then we may reduce this expression again for greater elegance,
with it being set $\frac{y}{1+y}=z$, so that it may thus be
$y=\frac{z}{1-z}$. Then presently the left member,
on account of $1+2y=\frac{1+z}{1-z}$, will admit this form:
$\frac{nz[(1+z)^n+(1-z)^n]}{(1+z)^n-(1-z)^n}$, because of which
it may therefore be equated to this continued fraction:
\[
  1 + \frac{(nn-1)zz}{3
          + \frac{(nn-4)zz}{5
          + \frac{(nn-9)zz}{7
          + \frac{(nn-16)zz}{9
          + \textrm{ etc.}}}}},
\]
which, by its elegance, merits the highest attention.

VII. Now, it is therefore manifest for this expression to always be halted
whenever $n$ is an integral number, either positive or negative.
It is moreover evident for the left member to retain the same value
even if for $n$ is written $-n$. Namely, this fact  
comes forth from:
\[
\frac{-nz[(1+z)^{-n}+(1-z)^{-n}]}{(1+z)^{-n}-(1-z)^{-n}},
\]
which fraction, if multiplied above and below by $(1-zz)^n$, induces this
form:
\[
\frac{-nz[(1-z)^n+(1+z)^n]}{(1-z)^n-(1+z)^n}=
\frac{nz[(1+z)^n+(1-z)^n]}{(1+z)^n-(1-z)^n},
\]
which is the same as the preceding expression. Thus it is the same
whether the positive or negative of the letter $n$ is taken.

VIII. Thus if we were to take $n=\pm 1$, the left member would be equal to 1,
which is moreover the value of the right. Furthermore, 
by putting $n=\pm 2$, the left member will come forth as equal to
$1+zz$, and indeed the right member will also be equal to
$1+zz$. In a similar way, by taking $n=\pm 3$, the left part,
and in turn the right, becomes $\frac{3(1+3zz)}{3+zz}$.

IX. Here one may deduce several conclusions of great importance,
depending on whenever a vanishing or infinite value is taken for the
exponent $n$, 
but first of all for the case in which an imaginary value is taken for the
letter $z$, which leads to an outstanding conclusion,
since this continued fraction shall nevertheless remain real,
for which conclusion we will therefore take up first.

\begin{center}
{\Large Conclusion I.\\
where $z=t\surd{-1}$}
\end{center} 

X. In this case therefore the continued fraction will have this form:
\[
  1 - \frac{(nn-1)tt}{3
          - \frac{(nn-4)tt}{5
          - \frac{(nn-9)tt}{7
          - \frac{(nn-16)tt}{9
          - \textrm{ etc.}}}}},
\]
and to be sure the left part will now be:
\[
\frac{nt\surd{-1}[(1+t\surd{-1})^n+(1-t\surd{-1})^n]}
{(1+t\surd{-1})^n-(1-t\surd{-1})^n},
\]
which not having been opposed by imaginary parts ought certainly
to have a real value, which we shall now investigate.
Then to this end we put $t=\frac{\sin{\phi}}{\cos{\phi}}$, so
that it will thus be $t=\tang{\phi}$; then it will therefore be:
\[
(1+t\surd{-1})^n=\frac{(\cos{\phi}+\surd{-1}\sin{\phi})^n}{\cos{\phi}^n}=
\frac{\cos{n\phi}+\surd{-1}\sin{n\phi}}{\cos{\phi}^n},
\]
and in a similar way:
\[
(1-t\surd{-1})^n=\frac{(\cos{\phi}-\surd{-1}\sin{\phi})^n}{\cos{\phi^n}}=
\frac{\cos{n\phi}-\surd{-1}\sin{n\phi}}{\cos{\phi}^n}.
\]
Therefore by substituting these values our left member will come forth:
\[
\frac{2n\surd{-1}\cdot \tg{\phi}\cos{n\phi}}{2\surd{-1}\sin{n\phi}}=
\frac{n\tg{\phi}\cos{n\phi}}{\sin{n\phi}}=\frac{n\tg{\phi}}{\tg{n\phi}}.
\]

XI. Therefore by putting $\tg \phi=t$ we will have the following
most remarkable continued fraction;
\[
  \frac{nt}{\tg n\phi} = 1 - \frac{(nn-1)tt}{3
          - \frac{(nn-4)tt}{5
          - \frac{(nn-9)tt}{7
          - \textrm{ etc.}}}},
\]
which then will be able to be represented in this way:
\[
  \tg n\phi= \frac{nt}{1
          - \frac{(nn-1)tt}{3
          - \frac{(nn-4)tt}{5
          - \frac{(nn-9)tt}{7 - \textrm{ etc.}}}}}, 
\] 
which 
expression therefore is able to be helpfully applied  to the tangents
of multiplied angles which are to be expressed by the tangent of
the single angle $t$. Thus if it were $n=2$, we will have
$\tg 2\phi=\frac{2t}{1-tt}$. In the very same way if $n=3$, it will be:
\[
\tg 3\phi=\frac{3t}{1-\frac{8tt}{3-tt}}=\frac{3t-t^3}{1-3tt}.
\]
In the most notable case presented whenever the exponent $n$ is taken as
less than infinite, 
it will then be $\tg n\phi=n\phi$,
therefore, by dividing both sides by $n$, this form arises:
\[
  \phi = \frac{t}{1
          + \frac{tt}{3
          + \frac{4tt}{5
          + \frac{9tt}{7
          +  \textrm{ etc.}}}}}, 
\] 
where the continued fraction is expressed by the tangent $t$ of the
angle itself.

XII. We will now consider the case in which an infinite magnitude is taken
for the exponent $n$, but when the angle $\phi$ less than infinite,
and then too for the tangent $t$ of it less than infinite, so that it would
thus be $n\phi=\theta$, and then also $nt=\theta$; then we will
therefore have such a continued fraction:
\[
  \tg \theta = \frac{\theta}{\displaystyle 1
          - \frac{\theta \theta}{\displaystyle 3
          - \frac{\theta \theta}{\displaystyle 5
          - \frac{\theta \theta}{\displaystyle 7
          - \textrm{ etc.}}}}}, 
\] 
by which formula, from the given angle $\theta$,
the tangent of it will be able to be determined,
which expression will be able to be seen just as the reciprocal of the
preceding.

\begin{center}
{\Large Conclusion II.\\
where a vanishing exponent $n$ is taken:}
\end{center}

XIII. Therefore in this case the continued fraction will be:
\[
  1 - \frac{zz}{\displaystyle 3
          - \frac{4zz}{\displaystyle 5
          - \frac{9zz}{\displaystyle 7
          - \frac{16zz}{\displaystyle 9
          - \textrm{ etc.}}}}}.
\] 
It is moreover to be noted for the left part to be 
$\frac{(1+z)^n-1}{n}=l(1+z)$, and for that reason
$(1+z)^n=1+nl(1+z)$; in a similar way it will be:
$(1-z)^n=1+nl(1-z)$, from which the left member comes forth as
\[
nz\frac{[2+nl(1+z)+nl(1-z)]}{nl(1+z)-nl(1-z)}=\frac{2z}{l\frac{1+z}{1-z}};
\]
here therefore we will have such a form:
\[
  \frac{2z}{l\frac{1+z}{1-z}}=1 - \frac{zz}{\displaystyle 3
          - \frac{4zz}{\displaystyle 5
          - \frac{9zz}{\displaystyle 7
          - \frac{16zz}{\displaystyle 9
          - \textrm{ etc.}}}}}, 
\] 
and then this logarithm may be expressed in the following way:
\[
  l\frac{1+z}{1-z} = \frac{2z}{\displaystyle 1
          - \frac{zz}{\displaystyle 3
          - \frac{4zz}{\displaystyle 5
          - \textrm{ etc.}}}}
\] 

\begin{center}
{\Large Conclusion III.\\
where an infinite magnitude is taken for the exponent $n$}
\end{center}

XIV. Here therefore, so that a finite value may be obtained for this
continued fraction, of which having been advanced a quantity less than infinity is to be taken for $z$,
it is put $nz=v$, so that it would be $z=\frac{v}{n}$, and thus our
continued fraction will be:
\[
  1 + \frac{vv}{3
          + \frac{vv}{5
          + \frac{vv}{7
          + \frac{vv}{9 + \textrm{ etc.}}}}}.
\]
Also, it is apparent for the left member to be $(1+\frac{v}{n})^n=e^v$,
and in a similar way $(1-\frac{v}{n})^n=e^{-v}$; therefore the left
member will have this form:
\[
\frac{v(e^v+e^{-v}}{e^v-e^{-v}}=\frac{v(e^{2v}+1)}{e^{2v}-1},
\]
and from this fact we will have this remarkable continued fraction:
\[
  \frac{v(e^{2v}+1)}{e^{2v}-1} = 1 + \frac{vv}{3
          + \frac{vv}{5
          + \frac{vv}{7
          + \frac{vv}{9
          + \textrm{ etc.}}}}}, 
\] 
whose transcendent value is actually able to be exhibited in this way by
a series:

\[
\frac{1+\frac{vv}{1\cdot 2}+\frac{v^4}{1\cdot 2\cdot 3\cdot 4}+
\frac{v^6}{1\cdot 2\cdots 6}+\textrm{ etc.}}
{1+\frac{vv}{1\cdot 2\cdot 3}+\frac{v^4}{1\cdot 2\cdots 5}+
\frac{v^6}{1\cdot 2\cdots 7}+\textrm{ etc.}}
\]

\end{document}